\documentclass{elsarticle}

\usepackage[utf8]{inputenc}
\usepackage[greek,english]{babel}
\usepackage{alphabeta} 

\usepackage{graphicx}

\newcommand{\R}{\mathbb{R}}
\newcommand{\C}{\mathbb{C}}
\newcommand{\eat}[1]{}

\usepackage{enumitem}
\setlist{nolistsep,noitemsep}
\usepackage[hidelinks]{hyperref}
\usepackage{amsmath, nccmath}
\usepackage{amssymb}
\usepackage{mathabx}
\usepackage[english]{babel}
\usepackage{amsthm}
\usepackage{multicol}
\usepackage{comment}
\usepackage{tikz}
\usepackage{ifthen}
\usetikzlibrary{calc}
\DeclareMathOperator{\rank}{rank}
\usepackage{relsize}
\usepackage[ruled,vlined]{algorithm2e}
\newtheorem{theorem}{Theorem}[section]

\newtheorem{lemma}[theorem]{Lemma}
\newtheorem{proposition}[theorem]{Proposition}
\newtheorem{definition}[theorem]{Definition}
\newtheorem{remark}[theorem]{Remark}
\newtheorem{problem}[theorem]{Problem}
\newtheorem{example}[theorem]{Example}

\newcommand*{\Comb}[2]{{}^{#1}C_{#2}}%
\newcommand{\Chi}{\mbox{\put(0,3){\large $\chi$}\hspace*{2.4mm}}}
\newcommand{\norm}[1]{\left\lVert#1\right\rVert}

\newcommand{\Rn}{{\mathbb{R}^n}}
\newcommand{\Cn}{{\mathbb{C}^n}}

\newcommand{\Rnn}{\R^{n\times n}}

\newcommand{\Cadd}{{C_{\mbox{\footnotesize \rm add~}}}}
\newcommand{\Cresame}{{C_{\mbox{\footnotesize \rm re-same~}}}}
\newcommand{\Creincr}{{C_{\mbox{\footnotesize \rm re-incr~}}}}

\newcommand*{\Scale}[2][4]{\scalebox{#1}{$#2$}}

\journal{Journal}

\begin{document}
\begin{frontmatter}
\title{\LARGE \bf
Algebraic connectivity: local and global maximizer graphs}

\author[1]{Karim Shahbaz} 
\ead{karimshahee@ee.iitb.ac.in}
\author[1]{Madhu N Belur} 
\ead{belur@ee.iitb.ac.in}
\author[2]{Ajay Ganesh}
\ead{ganesh@ualberta.ca}

\address[1]{\small \rm Department of Electrical Engg, Indian Institute of Technology Bombay, India}
\address[2]{\small \rm Department of Chemical \& Material Engg, University of Alberta, Edmonton, Canada}

\begin{abstract}
Algebraic connectivity is one way to quantify graph connectivity, which in turn
gauges robustness as a network. 
In this paper, we consider the problem of
maximising algebraic connectivity both local and globally over all simple,
undirected, unweighted graphs
with a given number of vertices and edges. We pursue this optimization by equivalently minimizing the
largest eigenvalue of the Laplacian of the `complement graph'. 
We establish that the union of complete subgraphs are
largest eigenvalue \emph{local} minimizer graphs. Further, under sufficient conditions 
satisfied by the edge/vertex counts
we prove that this union of complete components graphs are, in fact,
Laplacian largest eigenvalue \emph{global} maximizers; these results generalize 
the ones in the literature that are for just two components. These sufficient conditions
can be viewed as quantifying situations where
the component sizes are either `quite homogeneous' or some of them are relatively `negligibly small',
and thus generalize known results of homogeneity of components.
We finally relate this optimization
with the Discrete Fourier Transform (DFT) and circulant graphs/matrices.
\end{abstract}

\end{frontmatter}

\section{Introduction}\label{sec:intro}

Graph connectivity finds application in networking, network security, transportation systems, multi-agent control
and has been well studied in the literature.
Connectivity of a graph $G$ is also a measure of robustness as a network.
Algebraic connectivity being one of measures of graph connectivity is defined as the
second smallest eigenvalue $\lambda_{n-1}$ of the Laplacian matrix
$L(G) \in \Rnn$ of the unweighted, undirected and simple graph $G$.
In this paper, we consider only simple undirected, unweighted graphs with no self loops
and no multiple edges between any pairs of the vertices.
We study the problem of maximizing the algebraic connectivity of a graph for a given number of nodes and edges.
We pursue this problem both: a {\em global} maximization across all graphs, and a {\em local} sense, in which
we consider only one edge `rearrangements'
(defined precisely in Definition~\ref{def:OneEdgeReconnect} below).
Since algebraic connectivity and the problem of maximizing has received extensive attention and is well-understood,
we quickly delve further into the problem formulation, and then touch other closely related work
in the literature.

\subsection{Notation} \label{subsec:notation}

The notation we follow is standard and is included here for quick reference. The sets of real and complex numbers
are denoted respectively by $\R$ and $\C$. The largest eigenvalue of a
symmetric matrix is denoted by $\lambda_1$.
Given an undirected graph $G$, the number of vertices $|V(G)|$ is usually $n$, the number of edges $|E(G)|$ is
usually $m$, and the number of components of the graph is usually $p$. 
Further, the maximum degree across all vertices is denoted by $\Delta$ and
$d_{avg}$ is the average degree of vertices. 
The $n$ eigenvalues of the Laplacian matrix $L(G)$ are denoted by 
$\lambda_1(L(G))\geqslant\lambda_2(L(G))\geqslant\cdots \geqslant
\lambda_{n-1}(L(G))\geqslant \lambda_n(L(G))=0$.
When the matrix $L(G)$ and the graph $G$ are clear from the context,
we use just $\lambda_1$, $\dots$, $\lambda_n$
to denote the eigenvalues, and when comparing the maximum eigenvalues
of Laplacian matrices
of different graphs, say $G^m$ and $G^c$, we use $\lambda_1(G^m)$ and $\lambda_1(G^c)$.
Note that, since $L$ is symmetric, $\lambda_1(L(G))=\underset{\norm{x}_{_{2}}=1}{\max} x^TL(G)x$.

We deal with integer-valued properties and their relation with other bounds, and in this
context, we use the standard floor and the ceiling functions of $x$, denoted by $\lfloor  x \rfloor$
and $\lceil x \rceil$, to mean the largest/smallest integer not greater-than/not-smaller-than
the real number $x$ respectively.

The complete graph in $n$ vertices is denoted by $K_n$, and the complete bipartite graph
with vertex sets having cardinalities $p$ and $q$ is denoted by $K_{p,q}$. Of course,
our paper deals with complete multi-partite graphs, and in fact, with their complement graphs:
which would then be union of complete graphs, denoted by $\bigcup K_i$.

The notion of complement graph $G^c$ of a graph $G^m$ is straightforward: it is a
simple undirected graph with the same number (and indexing) of nodes and
in which there is an edge
in $G^c$ between two nodes, by definition, if and only if there is no edge in $G^m$.

\subsection{Problem formulation}

The paper deals with the following mutually closely related problems.

\begin{problem} \label{prob:LELM:LEGM}
The following sub-problems are inter-related for reasons clarified soon in the next section.
\begin{enumerate}
\item[(a)] For a given number of vertices $|V|=n$ and 
number of edges $|E|=m_1$, find an algebraic connectivity maximiser graph $G_1=(V, E)$.
\item[(b)] For given number of vertices $|V|=n$ and number of edges $|E|=m_2$,
minimise the largest Laplacian eigenvalue of the graph $G_2=(V, E)$.
\end{enumerate}
Further, each of the above optimizations can be pursued in one of two ways: globally and locally.
For simplicity, we elaborate on just the second one, i.e. the largest eigenvalue minimization: we study 
the global case, and the `local' case. More precisely,
\begin{enumerate}
\item finding a Largest Eigenvalue \underline{Global} Minimizer (LEGM) graph 
that has the least largest eigenvalue possible for the given number of vertices and edges, and
\item finding a Largest Eigenvalue \underline{Local} Minimizer (LELM), with `local minima'
in the
sense that all one-edge reconnect graphs (see Definition~\ref{def:OneEdgeReconnect})
have either the same largest eigenvalue or higher.  
\end{enumerate}
\end{problem}

\noindent
Related work in the context of the above problem is pursued in the next section.
The problem we consider in this paper also has a close link with circulant graphs 
(pursued further in Section~\ref{sec:CirculantMatrices}) and
DFT of time-symmetric vectors with entries from $\{0,1\}$.
The remark after the problem formulation below makes this precise.

\begin{problem} \label{prob:DFT}
DFT magnitude minimization: Given positive
integers $d$ and $n$ with
$1\leqslant d \leqslant n-1 $, consider a
vector $x\in\{0,1\}^n$ with $x_1 = 0$ and $\| x \| _1 = d$, and further, $x$ being `time-symmetric',
i.e. $x_{i} = x_{n+2-i}$ for $i=2,\dots,n$. Define $\bar{x}\in\Rn$ using $x$ by
$\bar{x}_1 = -d$, and $\bar{x}_i = x_i$ for all other $i$. Define the Discrete Fourier
Transform (DFT)\footnote{For uniformity
   with the rest of this paper, we use indices of $x,\bar{x}\in\Rn$ and $X\in\Cn$ to vary from
   $1$ to $n$, notwithstanding the typical DFT convention of using indices from 0 to $n-1$ for
   $x,\bar{x}$ and $X$.}
of the vector $\bar{x}$ by $X=DFT(\bar{x})$, and notice that $X \in \Rn$ due to
the assumed time-symmetry.
Consider the minimization problem: find $x$
satisfying the conditions above such that $\| X \|_\infty$ is minimized.
\end{problem}

Circulant matrices are pursued further in Section~\ref{sec:CirculantMatrices}. The following
remark motivates the assumptions within the problem formulations above.

\begin{remark} \label{rem:DFT:Circulant}
The following points relate Problems~\ref{prob:LELM:LEGM}
and~\ref{prob:DFT} and Laplacian matrices of circulant graphs. 
.\begin{enumerate}
\item  The condition $\bar{x}_1 = -d$ means that the `DC part' of $\bar{x}$ is zero and hence
$X_1 = 0$. Thus
minimizing $\|X\|_\infty$ means that the focus is on the minimization of
the maximum magnitude of all frequencies, except the DC.
\item Entries in $X$ are nothing but the negative of the eigenvalues
of the Laplacian of the graph $G_C$ constructed from $x$, and $G_C$ is regular (of degree $d$) and
is circulant; i.e., the Laplacian matrix is a circulant matrix.
\item The operation of defining $\bar{x} \in \Rn$ from $x\in\{0,1\}^n$ is one of adding an appropriately
scaled discrete time impulse $\delta$; the impulse has equal amount of all frequencies.
The DFT operation being linear on the signal space, this thus keeps the optimization
focus on the non-DC part in the signal $x$.
\item The operation of defining $\bar{x}$ from $x$ is like studying the eigenvalues of
$A-D$ (i.e. $-L$) instead of the adjacency matrix $A$, and note that the
diagonal matrix $D$ (the degree matrix) is merely $d \cdot I$ for this regular and circulant graph.
\end{enumerate}
\end{remark}

\subsection{Organization of the paper}

The rest of this paper is organized as follows. 
The next section relates the problem we pursue with other work in the literature and
in what way our work generalizes existing results. Section~\ref{sec:LocallyOptGraphs}
contains the main results of this paper, about {\em locally} optimal graphs. Further,
in the context of {\em globally} optimal graphs, our main results 
that improve upon results in the literature and also formulate for the
case of many components are contained in Section~\ref{sec:GlobOptGraphs}.
In Section~\ref{sec:CirculantMatrices}, we relate our work to the Discrete Fourier
Transform and circulant matrices/graphs. We consider some examples
in Section~\ref{sec:examples}. We conclude the paper in Section~\ref{sec:conclusion}, where we also summarize the contribution in this paper.

\section{Background and other work in this area}\label{sec:intro1}

Algebraic connectivity maximization of graphs has received much attention.
The survey papers
\cite{MerrisI},
\cite{Merris}, \cite{MerrisII},
\cite{Nikiforov} and
\cite{UpperBound} contain a wealth
of results about upper/lower bounds on the algebraic connectivity, many of which
we use crucially in our paper too.
In particular, given that we 
pursue maximum eigenvalue minimization on the complement graph
instead of directly algebraic connectivity (second-smallest
eigenvalue) maximization, it would help the reader to quickly review Proposition~\ref{lem:1} below
to see why this approach of focussing on the complement is equivalent.  

In the context of weighted graphs, \cite{Kim} proposes
an algorithm to find an edge to add to the graph to maximise algebraic connectivity, however,
the edge weight here is a function of distance between the vertices. 
Closely aligned with our paper, \cite{Takahashi} pursues
both Algebraic Connectivity `Local' Maximizers (ACLM) in 
the graph set of all one edge changes as in Definition~\ref{def:OneEdgeReconnect} and also
global maximizers, where for a given number of vertices and edges, conditions
are formulated.  Propositions~\ref{lem:BipartiteACM} and \ref{lem:BipartiteACLM} contain
the exact statements from \cite{Takahashi}, since this work
is relevant to the main results in our paper.
Both local and global optima obtained in \cite{Takahashi}
pursue for the case when the complement has {\em two} components,
while our paper generalizes to the case when the complement has any number of components,
and also slightly improves the bounds for the case of two components.

Recall that for a graph $G=(V, E)$, with $V$ the vertex set and $E$ the edge set,
the Laplacian matrix is defined as $L(G)=D(G)-A(G)$ where $D(G)$ is the diagonal matrix
with diagonal entries being degree of vertices and $A(G)$ is the adjacency matrix of graph $G$. 
The second smallest eigenvalue of $L$ is defined as the 
algebraic connectivity of the graph $G$: see \cite{Fiedler}. This eigenvalue is also
called the Fiedler value.  The rest of this section contains results that we use and/or improve
upon in this paper.

The following result crucially relates eigenvalues of the Laplacian matrices of a graph $G^m$ 
and its complement $G^c$. 
\begin{proposition} \label{lem:1} 
(\cite[page~148]{MerrisI}) Let $G^m$ be a simple undirected, unweighted graph and $G^{c}$ be its complement. Then the largest eigenvalue of the graph $\lambda_{1}(G^m)$ satisfies, $\lambda_{1}(G^m) \leqslant n$. Further, the eigenvalues of the Laplacian matrices of $G^m$ and $G^c$ are related by  $\lambda_{i}(G^c)=n-\lambda_{n-i}(G^m)$ for i=1, ..., n-1 and $\lambda_n(G^m)=\lambda_n(G^c)=0$.
\end{proposition}

The next well-known result (from \cite{Luo}) gives a lower bound for the maximum eigenvalue and also
formulates the unique situation when the bound is tight.
\begin{proposition} \label{lem:2} 
(\cite[Theorem 3.19]{Luo}) Consider a connected graph $G$ with at~least one edge, vertex set $V(G)$ of  cardinality $n$.
Then the following hold. 
\begin{enumerate} 
\item[a)] The maximum eigenvalue of the Laplacian matrix of the graph satisfies $\lambda_{1}(L(G))\geqslant \Delta+1$.
\item[b)]   $\lambda_{1}(L(G)) = \Delta+1$ holds if and only if $ \Delta=n-1$, i.e., there exists a `star node' in $G$.
\end{enumerate}  
\end{proposition}
Of course, if a graph is not connected, then the above result can still be used by noting the 
obvious fact that the Laplacian matrix $L_F$ of the full graph is a block diagonal matrix composed
of that of the individual components, and hence the eigenvalues of $L_F$ are the union of the 
individual Laplacian matrices' eigenvalues.
The following result gives a different lower bound
for the maximum eigenvalue and also the situation when this bound is tight.

\begin{proposition}\label{lem:DominatingNum} (\cite[Theorem 3]{Nikiforov}) 
Let Graph $G$ with $n\geqslant2$ vertices and
domination\footnote{\label{foot:domination:number} The 
    domination number of a graph \gamma(G) is defined as the minimum 
    size of the subset of vertices which are adjacent to every other vertex 
    of the graph} number,
denoted by $\gamma$. Then, $\lambda_1(G)\geqslant\lfloor \frac{n}{\gamma}\rfloor$ and, further, 
equality 
holds if and only if $G=G_a\bigcup G_b$ such that: 
\begin{enumerate}
    \item $|G_a|=\lfloor \frac{n}{\gamma}\rfloor$ and $\gamma(G_a)=1$, and
    \item $\gamma(G_b)=\gamma(G)-1$ and $\lambda_1(G_b)\leqslant \lfloor \frac{n}{\gamma}\rfloor$.
\end{enumerate}
\end{proposition}

The main results in our paper generalize the following results from \cite{Takahashi} and we generalize
these results to the case of more than two components (in the complement graph). 
For a specified number of vertices and edges, \cite{Takahashi} studies
the problem of Algebraic Connectivity Maximizer (ACM) graph and local algebraic maximizer 
graphs. The precise statements are below.

\begin{proposition}  \label{lem:BipartiteACM} 
\cite[Theorem~3]{Takahashi}: For integers $a \in \mathbb{Z}^+$, if $a\leqslant \lceil\frac{n}{2}\rceil$ and $a-\frac{2a^2}{n}<1$, then for any $n\geqslant3$, the complete bipartite graph $K_{a, n-a}$ is ACM in graphs with $n$ vertices and $a(n-a)$ edges.
\end{proposition}

\begin{proposition}  \label{lem:BipartiteACLM} 
\cite[Theorem~6]{Takahashi}: For integers $a \in \mathbb{Z}^+$, if $a\leqslant \lceil\frac{n}{2}\rceil$, then the complete bipartite graph $K_{a, n-a}$ is ACLM in graphs with $n$ vertices and $a(n-a)$ edges.
\end{proposition}

\begin{proposition}\label{lem:completeComp} 
\cite[Theorem 3.1]{Zhang} Consider Graph $G$ with at~least one edge and 
independence\footnote{The independence 
     number of graph $\alpha(G)$ is defined as the cardinality of the 
     largest set of vertices of the graph with no edge connection 
     between them} number $\alpha(G)$.
Then, 
$\lambda_1(G)\geqslant \frac{n}{\alpha}$ and, further,  equality holds 
if and only if $\alpha$ is factor of $n$ and thus $G$ then 
has $\alpha$ components each being $K_{\frac{n}{\alpha}}$.
\end{proposition}

\noindent We prove in this paper that the complement graph $G^c$ made up of two complete components graph is LEGM under a very similar (and slightly relaxed) sufficient condition as compared to Proposition~\ref{lem:BipartiteACM}.
We also extend the result of complete two components to multi-components and prove that the
graph is LEGM under an appropriately generalized sufficient condition. This result (Theorem~\ref{thm:MultiCompLEGM} below) generalizes Proposition~\ref{lem:completeComp} in a certain sense. The notion of 
Algebraic Connectivity Local Maximizers (ACLM) graph was introduced in
\cite{Takahashi}.  The ACLM graph is the one in which if one edge is changed (i.e.
one edge is either removed or reconnected to a different set of vertices), then
its algebraic connectivity remains highest among all such `one edge changed' graphs.
ACLM graphs are thus not globally optimal,
but at least locally optimal topologies and hence also usually globally suboptimal.
In \cite{Takahashi}, it has been shown that the complete bipartite graph
$K_{a, n-a}$ is an Algebraic Connectivity Local Maximizers (ACLM) in $G$ for
$n$ vertices and $a(n-a)$ edges graphs
for $2\leqslant a \leqslant \lfloor \frac{n}{2} \rfloor$; we generalize this result for
the case that the complement graph has not just two components but in fact
any number of components.


\section{Main results: locally optimal graphs}\label{sec:LocallyOptGraphs}

In this section we present the main results of this paper which concern `locally' optimal graphs. The notion
of local is made precise in the definition below. This notion coincides with that of \cite{Takahashi}.
Local optimality is important when only simple rearrangements of the topology of a set of multi-agents,
for example, is allowed and complicated rearrangements are disallowed. It helps to at least
be locally optimal.  Of course, globally optimal configurations would also need to satisfy this, and
thus local optimality conditions are necessary conditions for global optimality too.

\begin{definition}\label{def:OneEdgeReconnect}
    \begin{enumerate}
        \item[(a)] \textbf{One edge reconnect of $G_0$}: Let $G_0(V, E_0)$ be a simple graph with $|V|=n$, and $|E_0|=m$. We define $G_1(V, E_1)$ be a \underline{one-edge reconnect} of $G_0$ if $G_1$ is also a simple graph and one or both of nodes of exactly one edge differ from that of $G_0$. Thus, we have one-edge reconnect if $G_1$ satisfies $|E_1|=m$ and $\rank(L_1-L_0)=2$.
        \item[(b)]\textbf{One edge addition}: By one edge addition, we mean adding an edge to a graph while keeping the graph simple.
    \end{enumerate}
 \end{definition}

\noindent
Using the above notion of one edge reconnects and one edge additions,
we define a local minimizer graph; this is w.r.t. the largest eigenvalue of the Laplacian.
 
\begin{definition}\label{def:LELM}
\underline{\textbf{Largest Eigenvalue Local Minimizer graph}}: A graph $G_0$ is called a
Largest Eigenvalue \underline{Local} Minimizer (LELM) graph if $G_0$ has the least value
of the Laplacian matrix's largest eigenvalue amongst all the simple graphs $G$
obtained from $G_0$ by either a one edge reconnect or a one edge addition.
\end{definition} 



In the context of various possibilities of an edge reconnection or addition, it helps to visualize the
case using a figure. We include various figures, and the proof techniques vary depending on these cases.
In summary: when we have a union of complete subgraphs, then, an extra edge or
an edge reconnection connects to complete components, and we make a distinction
about whether the maximum degree increases or remains same, and whether
the largest component (with vertex-size say $n_1$), or vertex-size
slightly smaller than the largest (of size $n_1-1$),
or further smaller was involved in the edge reconnection/addition.
This distinction is needed to prove the local minimality of
the graph proposed in Theorem~\ref{thm:MultiCompLELM}.

\begin{figure}[h!]
        \centering
        \begin{tikzpicture}[scale=0.55]
    \tikzstyle{every node}=[draw, shape=circle, scale=0.5, fill=black!10];
    \foreach \i in {1,...,5}
    {   \pgfmathsetmacro\startj{\i+1}
        \foreach \j in {\startj, ..., 6}
        {
            \draw[thick] ({0+1.5*cos(360.0/6.0*\i)},{0+1.5*sin(360.0/6.0*\i)})--({0+1.5*cos(360.0/6.0*\j)},{0+1.5*sin(360.0/6.0*\j)});
        }
    }
    \foreach \i in {1,...,6}
    {
        \node (v_\i) at ({0+1.5*cos(360.0/6.0*\i)},{0+1.5*sin(360.0/6.0*\i)}){\i};
    }
    
    \foreach \k in {7,...,10}
    {   \pgfmathsetmacro\startl{\k+1}
        \foreach \l in {\startl, ..., 10}
        {
            \draw[thick] ({2*2.5+1.5*cos(360.0/4.0*(\k-6))},{0+1.5*sin(360.0/4.0*(\k-6))})--({2*2.5+1.5*cos(360.0/4.0*(\l-6))},{0+1.5*sin(360.0/4.0*(\l-6))});
        }
    }
    \foreach \k in {7,...,10}
    {   
        \node (v_\k) at ({2*2.5+1.5*cos(360.0/4.0*(\k-6))},{0+1.5*sin(360.0/4.0*(\k-6)}){\k}; 
    }
    \foreach \k in {11,...,13}
    {   \pgfmathsetmacro\startl{\k+1}
        \foreach \l in {\startl, ..., 13}
        {
            \draw[thick] ({2*4.5+1.5*cos(360.0/3.0*(\k-10))},{0+1.5*sin(360.0/3.0*(\k-10))})--({2*4.5+1.5*cos(360.0/3.0*(\l-10))},{0+1.5*sin(360.0/3.0*(\l-10))});
        }
    }
    \foreach \k in {11,...,13}
    {
        \node (v_\k) at ({2*4.5+1.5*cos(360.0/3.0*(\k-10))},{0+1.5*sin(360.0/3.0*(\k-10)}){\k}; 
    }
    \draw[thick](v_10) -- (v_12);
    \tikzstyle{every node}=[scale=0.75, fill=black!0];
     \node ($G_1$) at (0,-3){$G_1$};
    \node ($G_i$) at (5,-3){$G_i$};
    \node ($G_j$) at (9,-3){$G_j$};
    \end{tikzpicture}
        \caption{Connection established by one edge addition between $G_i$ and $G_j$, where $|V(G_j)| \leqslant |V(G_i)|\leqslant |V(G_1)|-2$}
        \label{fig:Case1:1a}
    \end{figure}
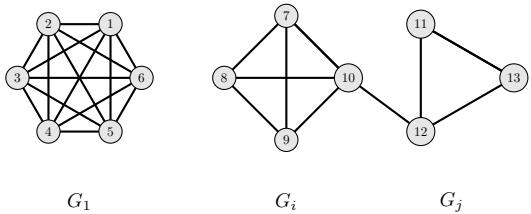
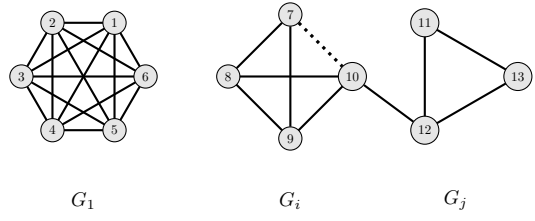
\begin{figure}[h!]
          \centering
          \begin{tikzpicture}[scale=0.55]
        \tikzstyle{every node}=[draw, shape=circle, scale=0.5, fill=black!10];
        \foreach \i in {1,...,5}
        {   \pgfmathsetmacro\startj{\i+1}
            \foreach \j in {\startj, ..., 6}
            {
                \draw[thick] ({0+1.5*cos(360.0/6.0*(\i))},{0+1.5*sin(360.0/6.0*(\i))})--({0+1.5*cos(360.0/6.0*(\j))},{0+1.5*sin(360.0/6.0*(\j))});
            }
        }
        \foreach \i in {1,...,6}
        {
            \node (v_\i) at ({0+1.5*cos(360.0/6.0*(\i))},{0+1.5*sin(360.0/6.0*(\i))}){\i};
        }
        
        \foreach \k in {8,...,9}
        {   \pgfmathsetmacro\startl{\k+1}
            \foreach \l in {\startl, ..., 10}
            {   
                    \draw[thick] ({2*2.5+1.5*cos(360.0/4.0*(\k-6))},{0+1.5*sin(360.0/4.0*(\k-6))})--({2*2.5+1.5*cos(360.0/4.0*(\l-6))},{0+1.5*sin(360.0/4.0*(\l-6))});
            }
        }
        \foreach \k in {7,...,10}
        {
            \node (v_\k) at ({2*2.5+1.5*cos(360.0/4.0*(\k-6))},{0+1.5*sin(360.0/4.0*(\k-6))}){\k};
        }
        \foreach \k in {11,12}
        {   \pgfmathsetmacro\startl{\k+1}
            \foreach \l in {\startl, ..., 13}
            {   
                \draw[thick] ({2*4.5+1.5*cos(360.0/3.0*(\k-10))},{0+1.5*sin(360.0/3.0*(\k-10))})--({2*4.5+1.5*cos(360.0/3.0*(\l-10))},{0+1.5*sin(360.0/3.0*(\l-10))});
            }
        }
        \foreach \k in {11,...,13}
        {
            \node (v_\k) at ({2*4.5+1.5*cos(360.0/3.0*(\k-10))},{0+1.5*sin(360.0/3.0*(\k-10))}){\k};
        }
        \draw[thick](v_10) -- (v_12);
        \draw[thick](v_7) -- (v_8);\draw[thick](v_7) -- (v_9);
        \draw[dotted, very thick] (v_7) -- (v_10);
        \tikzstyle{every node}=[scale=0.75, fill=black!0];
         \node ($G_1$) at (0,-3){$G_1$};
        \node ($G_i$) at (5,-3){$G_i$};
        \node ($G_j$) at (9,-3){$G_j$};
        \end{tikzpicture}
          \caption{Reconnection without increasing the maximum degree of $G_i$, where $|V(G_j)| \leqslant |V(G_i)|\leqslant |V(G_1)|-2$}
          \label{fig:Case1:1rs}
      \end{figure} 

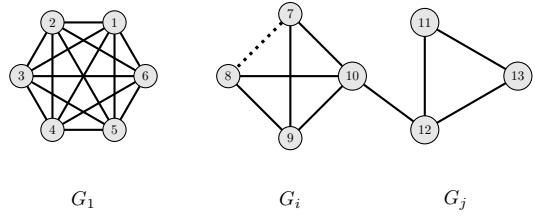
\begin{figure}[h!]
            \centering
            \begin{tikzpicture}[scale=0.55]
        \tikzstyle{every node}=[draw, shape=circle, scale=0.5, fill=black!10];
        \foreach \i in {1,...,5}
        {   \pgfmathsetmacro\startj{\i+1}
            \foreach \j in {\startj, ..., 6}
            {
                \draw[thick] ({0+1.5*cos(360.0/6.0*(\i))},{0+1.5*sin(360.0/6.0*(\i))})--({0+1.5*cos(360.0/6.0*(\j))},{0+1.5*sin(360.0/6.0*(\j))});
            }
        }
        \foreach \i in {1,...,6}
        {
            \node (v_\i) at ({0+1.5*cos(360.0/6.0*(\i))},{0+1.5*sin(360.0/6.0*(\i))}){\i};
        }
        
        \foreach \k in {8,...,9}
        {   \pgfmathsetmacro\startl{\k+1}
            \foreach \l in {\startl, ..., 10}
            {   
                    \draw[thick] ({2*2.5+1.5*cos(360.0/4.0*(\k-6))},{0+1.5*sin(360.0/4.0*(\k-6))})--({2*2.5+1.5*cos(360.0/4.0*(\l-6))},{0+1.5*sin(360.0/4.0*(\l-6))});
            }
        }
        \foreach \k in {7,...,10}
        {
            \node (v_\k) at ({2*2.5+1.5*cos(360.0/4.0*(\k-6))},{0+1.5*sin(360.0/4.0*(\k-6))}){\k};
        }
        \foreach \k in {11,12}
        {   \pgfmathsetmacro\startl{\k+1}
            \foreach \l in {\startl, ..., 13}
            {   
                \draw[thick] ({2*4.5+1.5*cos(360.0/3.0*(\k-10))},{0+1.5*sin(360.0/3.0*(\k-10))})--({2*4.5+1.5*cos(360.0/3.0*(\l-10))},{0+1.5*sin(360.0/3.0*(\l-10))});
            }
        }
        \foreach \k in {11,...,13}
        {
            \node (v_\k) at ({2*4.5+1.5*cos(360.0/3.0*(\k-10))},{0+1.5*sin(360.0/3.0*(\k-10))}){\k};
        }
        \draw[thick](v_10) -- (v_12);
        \draw[thick](v_7) -- (v_9);\draw[thick](v_7) -- (v_10);
        \draw[dotted, very thick] (v_7) -- (v_8);
        \tikzstyle{every node}=[scale=0.75, fill=black!0];
        \node ($G_1$) at (0,-3){$G_1$};
        \node ($G_i$) at (5,-3){$G_i$};
        \node ($G_j$) at (9,-3){$G_j$};
        \end{tikzpicture}
            \caption{Reconnection with increasing the maximum degree of $G_i$, where $|V(G_j)| \leqslant |V(G_i)|\leqslant |V(G_1)|-2$}
            \label{fig:Case1:1ri}
        \end{figure}
\begin{figure}[h!]
        \centering
        \begin{tikzpicture}[scale=0.55]
    \tikzstyle{every node}=[draw, shape=circle, scale=0.5, fill=black!10];
    \foreach \i in {1,...,5}
    {   \pgfmathsetmacro\startj{\i+1}
        \foreach \j in {\startj, ..., 6}
        {
            \draw[thick] ({0+1.5*cos(360.0/6.0*\i)},{0+1.5*sin(360.0/6.0*\i)})--({0+1.5*cos(360.0/6.0*\j)},{0+1.5*sin(360.0/6.0*\j)});
        }
    }
    \foreach \i in {1,...,6}
    {
        \node (v_\i) at ({0+1.5*cos(360.0/6.0*\i)},{0+1.5*sin(360.0/6.0*\i)}){\i};
    }
    
    \foreach \k in {7,...,10}
    {   \pgfmathsetmacro\startl{\k+1}
        \foreach \l in {\startl, ..., 11}
        {
            \draw[thick] ({2*2.5+1.5*cos(360.0/5.0*(\k-6))},{0+1.5*sin(360.0/5.0*(\k-6))})--({2*2.5+1.5*cos(360.0/5.0*(\l-6))},{0+1.5*sin(360.0/5.0*(\l-6))});
        }
    }
    \foreach \k in {7,...,11}
    {   
        \node (v_\k) at ({2*2.5+1.5*cos(360.0/5.0*(\k-6))},{0+1.5*sin(360.0/5.0*(\k-6)}){\k}; 
    }
    \foreach \k in {12,...,14}
    {   \pgfmathsetmacro\startl{\k+1}
        \foreach \l in {\startl, ..., 15}
        {
            \draw[thick] ({2*5+1.5*cos(360.0/4.0*(\k-11))},{0+1.5*sin(360.0/4.0*(\k-11))})--({2*5+1.5*cos(360.0/4.0*(\l-11))},{0+1.5*sin(360.0/4.0*(\l-11))});
        }
    }
    \foreach \k in {12,...,15}
    {
        \node (v_\k) at ({2*5+1.5*cos(360.0/4.0*(\k-11))},{0+1.5*sin(360.0/4.0*(\k-11)}){\k}; 
    }
    \draw[thick](v_11) -- (v_13);
    \tikzstyle{every node}=[scale=0.75, fill=black!0];
     \node ($G_1$) at (0,-3){$G_1$};
    \node ($G_i$) at (5,-3){$G_i$};
    \node ($G_j$) at (10,-3){$G_j$};
    \end{tikzpicture}
        \caption{Connection established by one edge addition to $G_i$, where $|V(G_j)| \leqslant |V(G_i)|=|V(G_1)|-1$}
        \label{fig:Case2:2a}
    \end{figure}
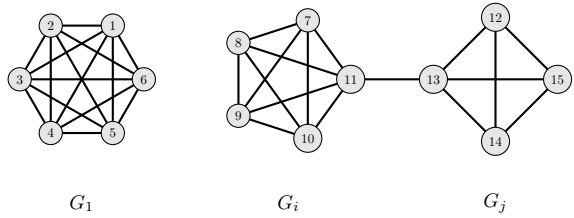
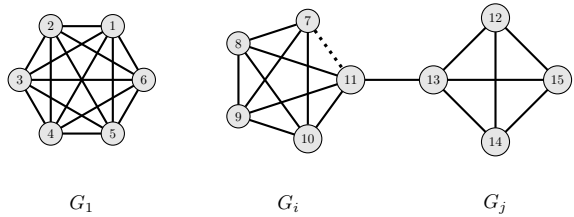
\begin{figure}[h!]
                 \centering
                 \begin{tikzpicture}[scale=0.55]
            \tikzstyle{every node}=[draw, shape=circle, scale=0.5, fill=black!10];
            \foreach \i in {1,...,5}
            {   \pgfmathsetmacro\startj{\i+1}
                \foreach \j in {\startj, ..., 6}
                {
                    \draw[thick] ({0+1.5*cos(360.0/6.0*\i)},{0+1.5*sin(360.0/6.0*\i)})--({0+1.5*cos(360.0/6.0*\j)},{0+1.5*sin(360.0/6.0*\j)});
                }
            }
            \foreach \i in {1,...,6}
            {
                \node (v_\i) at ({0+1.5*cos(360.0/6.0*\i)},{0+1.5*sin(360.0/6.0*\i)}){\i};
            }
            
            \foreach \k in {8,...,10}
            {   \pgfmathsetmacro\startl{\k+1}
                \foreach \l in {\startl, ..., 11}
                {
                    \draw[thick] ({2*2.5+1.5*cos(360.0/5.0*(\k-6))},{0+1.5*sin(360.0/5.0*(\k-6))})--({2*2.5+1.5*cos(360.0/5.0*(\l-6))},{0+1.5*sin(360.0/5.0*(\l-6))});
                }
            }
            \foreach \k in {7,...,11}
            {   
                \node (v_\k) at ({2*2.5+1.5*cos(360.0/5.0*(\k-6))},{0+1.5*sin(360.0/5.0*(\k-6)}){\k}; 
            }
            \foreach \k in {12,...,14}
            {   \pgfmathsetmacro\startl{\k+1}
                \foreach \l in {\startl, ..., 15}
                {
                    \draw[thick] ({2*5+1.5*cos(360.0/4.0*(\k-11))},{0+1.5*sin(360.0/4.0*(\k-11))})--({2*5+1.5*cos(360.0/4.0*(\l-11))},{0+1.5*sin(360.0/4.0*(\l-11))});
                }
            }
            \foreach \k in {12,...,15}
            {
                \node (v_\k) at ({2*5+1.5*cos(360.0/4.0*(\k-11))},{0+1.5*sin(360.0/4.0*(\k-11)}){\k}; 
            }
            \draw[thick](v_7) -- (v_8); \draw[thick](v_7) -- (v_9); \draw[thick](v_7) -- (v_10);
            \draw[dotted, very thick] (v_7) -- (v_11); \draw[thick](v_11) -- (v_13);
            \tikzstyle{every node}=[scale=0.75, fill=black!0];
            \node ($G_1$) at (0,-3){$G_1$};
            \node ($G_i$) at (5,-3){$G_i$};
            \node ($G_j$) at (10,-3){$G_j$};
            \end{tikzpicture}
                 \caption{Reconnection without increasing the maximum degree of $G_i$, where $|V(G_j)| \leqslant |V(G_i)|=|V(G_1)|-1$}
                 \label{fig:Case2:2rs}
             \end{figure}
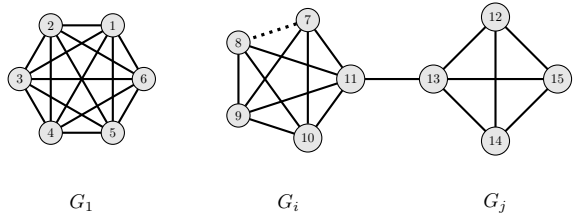
\begin{figure}[h!]
            \centering
            \begin{tikzpicture}[scale=0.55]
            \tikzstyle{every node}=[draw, shape=circle, scale=0.5, fill=black!10];
            \foreach \i in {1,...,5}
            {   \pgfmathsetmacro\startj{\i+1}
                \foreach \j in {\startj, ..., 6}
                {
                    \draw[thick] ({0+1.5*cos(360.0/6.0*\i)},{0+1.5*sin(360.0/6.0*\i)})--({0+1.5*cos(360.0/6.0*\j)},{0+1.5*sin(360.0/6.0*\j)});
                }
            }
            \foreach \i in {1,...,6}
            {
                \node (v_\i) at ({0+1.5*cos(360.0/6.0*\i)},{0+1.5*sin(360.0/6.0*\i)}){\i};
            }
            
            \foreach \k in {8,...,10}
            {   \pgfmathsetmacro\startl{\k+1}
                \foreach \l in {\startl, ..., 11}
                {
                    \draw[thick] ({2*2.5+1.5*cos(360.0/5.0*(\k-6))},{0+1.5*sin(360.0/5.0*(\k-6))})--({2*2.5+1.5*cos(360.0/5.0*(\l-6))},{0+1.5*sin(360.0/5.0*(\l-6))});
                }
            }
            \foreach \k in {7,...,11}
            {   
                \node (v_\k) at ({2*2.5+1.5*cos(360.0/5.0*(\k-6))},{0+1.5*sin(360.0/5.0*(\k-6)}){\k};
            }
            \foreach \k in {12,...,14}
            {   \pgfmathsetmacro\startl{\k+1}
                \foreach \l in {\startl, ..., 15}
                {
                    \draw[thick] ({2*5+1.5*cos(360.0/4.0*(\k-11))},{0+1.5*sin(360.0/4.0*(\k-11))})--({2*5+1.5*cos(360.0/4.0*(\l-11))},{0+1.5*sin(360.0/4.0*(\l-11))});
                }
            }
            \foreach \k in {12,...,15}
            {
                \node (v_\k) at ({2*5+1.5*cos(360.0/4.0*(\k-11))},{0+1.5*sin(360.0/4.0*(\k-11)}){\k};
            }
            \draw[dotted, very thick] (v_7) -- (v_8); \draw[thick](v_7) -- (v_9); \draw[thick](v_7) -- (v_10);
            \draw[thick] (v_7) -- (v_11); \draw[thick](v_11) -- (v_13);
            \tikzstyle{every node}=[scale=0.75, fill=black!0];
            \node ($G_1$) at (0,-3){$G_1$};
            \node ($G_i$) at (5,-3){$G_i$};
            \node ($G_j$) at (10,-3){$G_j$};
            \end{tikzpicture}
            \caption{Reconnection with increasing maximum degree of $G_i$, where $|V(G_j)| \leqslant |V(G_i)|=|V(G_1)|-1$}
            \label{fig:Case2:2ri}
        \end{figure}
\begin{figure}[h!]
        \centering
        \begin{tikzpicture}[scale=0.45]
    \tikzstyle{every node}=[draw, shape=circle, scale=0.5, fill=black!10];
    
    \foreach \i in {1,...,5}
    {\pgfmathsetmacro\startj{\i+1}
        \foreach \j in {\startj, ..., 6}
        {
            \draw[thick] ({0+2*cos(360.0/6.0*(\i))},{0+2*sin(360.0/6.0*(\i)})--({0+2*cos(360.0/6.0*(\j))},{0+2*sin(360.0/6.0*(\j)});
        }
    }
    \foreach \i in {1,...,6}
    {
        \node (v_\i) at ({0+2*cos(360.0/6.0*(\i))},{0+2*sin(360.0/6.0*(\i)}){\i};
    }
    
    \foreach \k in {7,...,11}
    {\pgfmathsetmacro\startl{\k+1}
        \foreach \l in {\startl, ..., 12}
        {   \pgfmathsetmacro\coordk{\k-6}
            \pgfmathsetmacro\coordl{\l-6}
            \draw[thick] ({2*3+2*cos(360.0/6.0*\coordk)},{0+2*sin(360.0/6.0*\coordk)}) -- ({2*3+2*cos(360.0/6.0*\coordl)},{0+2*sin(360.0/6.0*\coordl)});
        }
    }
    \foreach \k in {7,...,12}
    {
        \node (v_\k) at ({2*3+2*cos(360.0/6.0*(\k-6))},{0+2*sin(360.0/6.0*(\k-6)}){\k};
    }
    \draw[thick] (v_6) -- (v_9);
    \tikzstyle{every node}=[scale=0.75, fill=black!0];
     \node ($G_1$) at (0,-3){$G_1$};
    \node ($G_i$) at (6,-3){$G_i$};
    \end{tikzpicture}
        \caption{Connection established between largest size component $G_1$ with any other size component by one edge addition}
        \label{fig:Case3:3a}
    \end{figure}
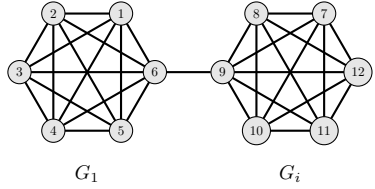
\begin{figure}[h!]
    \centering
    \begin{tikzpicture}[scale=0.45]
\tikzstyle{every node}=[draw, shape=circle, scale=0.5, fill=black!10];

\foreach \i in {2,...,5}
{\pgfmathsetmacro\startj{\i+1}
    \foreach \j in {\startj, ..., 6}
    {
        \draw[thick] ({0+2*cos(360.0/6.0*(\i))},{0+2*sin(360.0/6.0*(\i)})--({0+2*cos(360.0/6.0*(\j))},{0+2*sin(360.0/6.0*(\j)});
    }
}
\foreach \i in {1,...,6}
{
    \node (v_\i) at ({0+2*cos(360.0/6.0*(\i))},{0+2*sin(360.0/6.0*(\i)}){\i};
}

\foreach \k in {7,...,10}
{\pgfmathsetmacro\startl{\k+1}
    \foreach \l in {\startl, ..., 11}
    {   \pgfmathsetmacro\coordk{\k-6}
        \pgfmathsetmacro\coordl{\l-6}
        \draw[thick] ({2*3+2*cos(360.0/5.0*\coordk)},{0+2*sin(360.0/5.0*\coordk)}) -- ({2*3+2*cos(360.0/5.0*\coordl)},{0+2*sin(360.0/5.0*\coordl)});
    }
}
\foreach \k in {7,...,11}
{
    \node (v_\k) at ({2*3+2*cos(360.0/5.0*(\k-6))},{0+2*sin(360.0/5.0*(\k-6)}){\k};
}
\draw[thick] (v_1) -- (v_2); \draw[thick] (v_1) -- (v_3); \draw[thick] (v_1) -- (v_4); 
\draw[thick] (v_1) -- (v_5); \draw[dotted, very thick] (v_1) -- (v_6); \draw[thick] (v_6) -- (v_9);
\tikzstyle{every node}=[scale=0.75, fill=black!0];
     \node ($G_1$) at (0,-3){$G_1$};
    \node ($G_i$) at (6,-3){$G_i$};
\end{tikzpicture}
    \caption{Reconnection of edge without increasing the maximum degree of largest size component $G_1$}
    \label{fig:Case3:3r}
\end{figure}
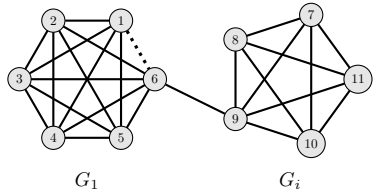

\begin{lemma}\label{lem:AddEdge}
Suppose a connection is established between complete graph components $G_i$ and $G_j$ by \underline{adding} an edge to give $G^{+}_{ij}$ and let $L_{new}-L_{old}=: \Cadd $ is connection matrix. Then  $\rank(L_{new}-L_{old})=1$ and the largest eigenvalue of $ \Cadd $ is $2$ (refer to Figure \ref{fig:Case1:1a}).
\end{lemma}
\begin{proof}
Contribution to the Laplacian matrix of graph
due to an edge addition has the structure: 
\begin{flushleft}
$ \Cadd =\left[ \begin{array}{c c | c c }
  &\textbf{$G_i$}~~~   &~~~\textbf{$G_j$}&\\ 
0 &0 &0 &0 \\
0 &1 &-1 &0\\
\cline{1-4}
0 &-1 &1 &0\\
 0 &0 &0 &0 \\
\end{array}\right]$.
\end{flushleft}
Clearly, the matrix $ \Cadd $ has $\rank$ one
and the characteristic polynomial:\\
\mbox{$\Chi_{ \Cadd }(s)=s^3(s-2)$}. So, $\lambda_1( \Cadd )=2$.\\
For bigger or general size $G_i$ and $G_j$ with $|V(G_i)|+|V(G_j)|=a$, the structure of $ \Cadd $ remains same but with zeros padded appropriately. Thus, $ \Cadd $ has $\rank$ one in general also and the lemma is proved. 
\end{proof}

\begin{lemma}\label{lem:ReconnectEdge_Idmax}
Suppose a connection is established between complete graph components $G_i$ and $G_j$ by \underline{reconnecting} an edge by removing one edge $e$ and adding elsewhere such that both nodes of $e$ change, to give $G^{+}_{ij}$ and let $L_{new}-L_{old}=:\Creincr$ is connection matrix.
Then  $\rank(L_{new}-L_{old})=2$ and the largest eigenvalue of $\Creincr$ is $2$ (refer to Figure \ref{fig:Case1:1ri}).
\end{lemma}
\begin{proof}
Contribution to the Laplacian matrix of graph due to
an edge reconnection as specified in the lemma has the following structure:
\begin{flushleft}
$\Creincr=\left[ \begin{array}{c c c | c}
 &  \textbf{$G_i$}& & \textbf{$G_j$} \\ 
-1 &1 &0 &0\\
1 &-1 &0 &0\\
0 &0 &1 &-1 \\
\cline{1-4}
0 &0 &-1 &1\\
\end{array}\right]$.
\end{flushleft}
Clearly, the matrix $\Creincr$ has $\rank$ two and the characteristic polynomial:\\
\mbox{$\Chi_{\Creincr}(s)=s^2(s^2-2^2)$}. Thus, $\lambda_1(\Creincr )=2$.
Again, for the general case, zeros get padded appropriately and the lemma is thus proved.
\end{proof}

\begin{lemma}\label{lem:ReconnectEdge_Sdmax}
Suppose a connection is established between complete graph components $G_i$ and $G_j$ by reconnecting an edge by removing one edge $e$ and adding an edge such that only one node of $e$ gets change, to give $G^{+}_{ij}$, and let $L_{new}-L_{old}=:\Cresame$ be the connection matrix. Then $\rank(L_{new}-L_{old})=2$ and the largest eigenvalue of $ \Cresame $ is $\sqrt{3}$ (refer to Figure \ref{fig:Case1:1rs}).
\end{lemma}
\begin{proof} 
Contribution to the Laplacian matrix of graph due to
an edge reconnection as specified in the lemma has the following structure: 
\begin{flushleft}
    $ \Cresame  = \left[ \begin{array}{c c c | c c}
 &  \textbf{$G_i$}  & & \textbf{$G_j$} &\\ 
 0 &0 &0 &0 &0\\
 0 &-1 &1 &0 &0\\
 0 &1 &0 &-1 &0\\
\cline{1-5}
 0 &0 &-1 &1 &0\\
 0 &0 &0 &0 &0\\
\end{array}\right]$. 
\end{flushleft}
Clearly, the matrix $ \Cresame $ has $\rank$ two and the
characteristic polynomial: \\
$\Chi_{ \Cresame }(s)=s^3(s^2-3)$.
Thus, $\lambda_1( \Cresame )=\sqrt{3}$.
Again, for the general case, zeros get padded appropriately and the lemma is thus proved.
\end{proof}
With the above lemmas, we are ready to state and prove the first main result of this paper. 

\begin{theorem}\label{thm:MultiCompLELM}: A graph $G=(V,E)$ which is a union of complete components is a Largest Eigenvalue \underline{Local} Minimizer (LELM).
In other-words, for graph $G$ of $n$ number of vertices, $m$ number of edges and $p$ number of
complete components of $|V(G_{i})|$ size such that $\sum_{i=1}^p |V(G_i)|=n$
and $m=\sum_{i=1}^p\Comb{|V(G_i)|}{2}$, then $\lambda_{1}(G)$ is locally minimized,
i.e. minimized w.r.t. one edge reconnects, one edge removals and
one edge additions (as defined in Definition \ref{def:OneEdgeReconnect}).
Further, $\lambda_1(G)=\underset{i\in \{1, 2,\cdots, p\}}{max}\{|V(G_i)|\}$.
\end{theorem}
\begin{proof}
Let $G_1, G_2, ..., G_p$ be the components and the number of nodes involved in those components
be $|V(G_i)|=n_i$ then $|E(G)|=\sum_{i\in \{1,2,..,p\}}\Comb{|V(G_i)|}{2}$. 
Without loss of generality, we assume that the components size have the
following relation between them: $|V(G_1)|\geqslant |V(G_2)|\geqslant |V(G_3)|\geqslant ... |V(G_p)| > 0$.
Thus the largest eigenvalue $\lambda_{1}(G)$ of the graph $G$ is: 
$\lambda_{1}\left(G\right)=|V(G_1)|$, since $\lambda_1(L(K_{n_1}))=n_1$.
In this setup if one edge is reconnected or one edge is added, it can be
connected in the following 3 ways: \\
\textbf{Case~1}: Between components of smaller sizes $G_i$ and $G_j$ such
that $|V(G_j)|\leqslant |V(G_i)|\leqslant |V(G_1)|-2$, i.e. both components $G_i$ and $G_j$
are at~least two or more nodes smaller than the largest component's size ($G_1$).\\
\noindent \textbf{Case~2}: Between component $G_i$ and $G_j$ with $|V(G_j)|\leqslant |V(G_i)|=|V(G_1)|-1$.\\
\noindent \textbf{Case~3}: Between $G_1$ and any other component: same size as $G_1$ or smaller.\\
\noindent We now prove the theorem for each of the 3 cases. Note that for each case, we have three subcases:
($a$) Addition of an edge, ($r_s$) Removal and addition of an edge $e$ such that only one vertex of $e$ is changed, and ($r_i$) Removal and addition of an edge $e$ such that both vertices of $e$ are changed.
We are not mentioning the one edge removal for local minimizer explicitly because by removing only one edge from complete component graphs, does not change its Laplacian largest eigenvalue.
Hence, this proposed graph is trivially LELM.\\
\noindent \textbf{Case~1}: Between components of smaller sizes $G_i$ and $G_j$ (without loss of generality assuming $|V(G_i)| \geqslant |V(G_j)|$) such that $|V(G_i)|\leqslant |V(G_1)|-2$, i.e. both $G_i$ \& $G_j$ are at~least two nodes smaller than the largest component $G_1$:
\begin{enumerate}
    \item[$1a$)] By one edge addition (refer to Figure \ref{fig:Case1:1a}): If connection is established between $G_i$ and $G_j$ to give $G^{+}_{ij}$ by adding an edge, then the connection matrix $ \Cadd $ of Lemma \ref{lem:AddEdge}, gets added to $L(G_i \oplus G_j)$.\\
Thus, due to the edge addition in between components we get,  $L(G^{+}_{ij})=L(G_i \oplus G_j)+ \Cadd $.\\
Also, $\lambda_1(L(G^{+}_{ij}))=\underset{\norm{x}_{_{2}}=1}{\max}x^TL(G^{+}_{ij})x=\underset{\norm{x}_{_{2}}=1}{\max}[x^TL(G_i\oplus G_j)x+x^T \Cadd x]$.
    
Using Lemma \ref{lem:AddEdge}, we have $\lambda_1( \Cadd )=2$, which 
implies that $\lambda_1(L(G^{+}_{ij})) \leqslant \lambda_1(L(G_i \oplus G_j))+ 2=\lambda_1(L(G_i))+ 2 \leqslant \lambda_1(G_1)$.
    
\noindent 
Therefore, $\lambda_1(G)=\max\{\lambda_1(G_1), \lambda_1(G_2), ..., \lambda_1(G^{+}_{ij})\}=\lambda_1(G_{1})$.  This proves that $\lambda_1(G)$ remains same
and the proposed graph $G$ is a $\lambda_{1}(G)$ local minimizer. 
    \item[$1r$)] One edge reconnect: If the connection established between $G_i$ and $G_j$ to give $G^{+}_{ij}$ by reconnecting one edge, then the following two different types of $C$ connection matrix get added to $L(G_i \oplus G_j)$ depending upon how the reconnection of edge is done.
        \begin{enumerate}
            \item[$1r_s$)] Reconnection \underline{without increasing the maximum degree} of $G_i$ (refer to Figure \ref{fig:Case1:1rs}):\\
            Due to the reconnection, the connection matrix $\Cresame$ of Lemma~\ref{lem:ReconnectEdge_Sdmax} gets added and we get $L(G^{+}_{ij})=L(G_i \oplus G_j)+\Cresame$.\\
        Also, $\lambda_1(L(G^{+}_{ij}))=\underset{\norm{x}_{_{2}}=1}{\max}x^TL(G^{+}_{ij})x=\underset{\norm{x}_{_{2}}=1}{\max}[x^TL(G_i\oplus G_j)x+x^T \Cresame x]$
        
        Using Lemma \ref{lem:ReconnectEdge_Sdmax}, we have $\lambda_1( \Cresame )=\sqrt{3}$.\\
        $\implies \lambda_1(L(G^{+}_{ij})) \leqslant \lambda_1(L(G_i \oplus G_j))+ \sqrt{3}=\lambda_1(L(G_i))+ \sqrt{3}<\lambda_1(L(G_i))+ 2 \leqslant \lambda_1(G_1)$.
        
        \noindent Therefore, $\lambda_1(G^{+})=\max\{\lambda_1(G_1), \lambda_1(G_2), ..., \lambda_1(G^{+}_{ij})\}=\lambda_1(G_{1})$.\\
        $\lambda_1(G)$ remains same and our graph is local minimizer.
      
        \item[$1r_i$)] Reconnection with \underline{increasing the maximum degree} of $G_i$ (refer to Figure \ref{fig:Case1:1ri}): \\
        Due to reconnection, the connection matrix $ \Creincr $ of Lemma~\ref{lem:ReconnectEdge_Idmax} gets added and we get  $L(G^{+}_{ij})=L(G_i \oplus G_j)+ \Creincr $.\\
        Also, $\lambda_1(L(G^{+}_{ij}))=~\underset{\norm{x}_{_{2}}=1}{\max}x^TL(G^{+}_{ij})x  =\underset{\norm{x}_{_{2}}=1}{\max}[x^TL(G_i\oplus G_j)x+x^T \Creincr x]$
        
        Using Lemma \ref{lem:ReconnectEdge_Idmax}, $\lambda_1( \Creincr )=2$.\\
        $\implies \lambda_1(L(G^{+}_{ij})) \leqslant \lambda_1(L(G_i \oplus G_j))+ 2=\lambda_1(L(G_i))+ 2 \leqslant \lambda_1(G_1)$.
        
        \noindent Therefore, $\lambda_1(G^+)=\max\{\lambda_1(G_1), \lambda_1(G_2), ..., \lambda_1(G^{+}_{ij})\}=\lambda_1(G_{1})$.\\
        $\lambda_1(G)$ remains same and the proposed graph $G$ graph is local minimizer. This completes the proof of Case~1.
        \end{enumerate}
\end{enumerate}
\noindent \textbf{Case~2}: Between component $G_i$ of size $|V(G_i)|=|V(G_1)|-1$ and any
other component $G_j$ of equal or smaller size than $G_i$ i.e. $|V(G_j)|\leqslant|V(G_i)|=|V(G_1)|-1$:
\begin{enumerate} 
    \item[$2a$)] By one edge addition (refer to Figure \ref{fig:Case2:2a}): Suppose connection is established between $G_i$ and $G_j$ to give $G^{+}_{ij}$ by adding an edge (using Proposition \ref{lem:2}~b),\\
    $\lambda_1(G^{+}_{ij}) > |V(G_i)|+1=|V(G_1)|=\lambda_1(G_1)$.\\
    $\lambda_1(G^{+})=\max\{\lambda_1(G_1), \lambda_1(G_2), ..., \lambda_1(G^{+}_{ij})\}=\lambda_1(G^{+}_{ij})>\lambda_1(G_1)$. Thus, proposed graph $G$ is a local minimizer. 
    \item[$2r$)] One edge reconnect: Suppose connection is established between $G_i$ and $G_j$ to give $G^{+}_{ij}$ by relocating an edge, then following two different type of $C$  connection matrix gets added to $L(G_i \oplus G_j)$ depending upon how reconnection of edge is done.
    \begin{enumerate} 
        \item[$2r_s$)] Reconnection \underline{without increasing the maximum degree} of $G_i$ (refer to Figure \ref{fig:Case2:2rs}):\\
            Due to reconnection, we get $L(G^{+}_{ij})=L(G_i \oplus G_j)+ \Cresame $.\\ 
            Also, $\lambda_1(L(G^{+}_{ij}))=\underset{\norm{x}_{_{2}}=1}{\max}x^TL(G^{+}_{ij})x=\underset{\norm{x}_{_{2}}=1}{\max}[x^TL(G_i\oplus G_j)x+x^T \Cresame x]$\\
            Using Lemma \ref{lem:ReconnectEdge_Sdmax},  $\lambda_1( \Cresame )=\sqrt{3}$.\\
            $\lambda_1(L(G^{+}_{ij})) \leqslant \lambda_1(L(G_i \oplus G_j))+ \sqrt{3}=\lambda_1(L(G_i))+\sqrt{3}=\lambda_1(L(G_1))+ \sqrt{3}-1$.\\
            \noindent So, in case of reconnecting without increasing maximum degree, we use the following relation:\\
$\lambda_1(G_i)=\lambda_1(G_1)-1<\lambda_1(G^{+}_{ij}) \leqslant \lambda_1(G_1)+\sqrt{3}-1$.\\
Thus, $\lambda_1(G_1)\leqslant \lambda_1(G^+) \leqslant \lambda_1(G_1)+\sqrt{3}-1$ implies $\lambda_1(G^+)$ either increases or remains same. Therefore again the proposed graph $G$ graph is an LELM.
            
\item[$2r_i$)] Reconnection \underline{with increasing the maximum degree} of $G_i$ (refer to Figure \ref{fig:Case2:2ri}): Suppose connection is established between $G_i$ and $G_j$ to give $G^{+}_{ij}$ by reconnecting an edge with increasing maximum degree of $G_i$, we get: \\(using Proposition \ref{lem:2}),  $\lambda_1(G^{+}_{ij}) > |V(G_i)|+1=|V(G_1)|=\lambda_1(G_1)$.\\
        $\lambda_1(G^{+})=\max\{\lambda_1(G_1), \lambda_1(G_2), ..., \lambda_1(G^{+}_{ij})\}=\lambda_1(G^{+}_{ij})$. Hence, $\lambda_1(G)$ increases. Thus, the proposed graph $G$ is an LELM.
\end{enumerate}
This completes the proof of Case 2.
\end{enumerate}

\noindent \textbf{Cases~3}: Between $G_1$ and any other component: 
\begin{enumerate}
    \item[$3a$)] By one edge addition (refer to Figure \ref{fig:Case3:3a}): Before addition of edge, we have $\lambda_{1}(G)=|V(G_1)|$. Then connection is established in two ways: between two largest size components and between largest and any other size components. Thus, addition of edge between components  $K_{|V(G_1)|}$ and $K_{|V(G_i)|}$ leads to $\lambda_{1} (G) > |V(G_1)|$ (using Proposition \ref{lem:2})\cite{Luo}, \cite{MerrisII}). Therefore, the proposed graph is a $\lambda_{1}(G)$ local minimizer (LELM).
    \item[$3r$)] Reconnecting of edge with or without increasing the maximum degree of $G_1$ (refer to Figure \ref{fig:Case3:3r}): Here, the connection is established in cases with largest size component $G_1$ by re-connecting $K_{|V(G_1)|}$ and $K_{|V(G_i)|}$ either by increasing maximum degree of $G_1$ or not; similarly like addition of edge, the reconnection leads to $\lambda_{1} (G) > |V(G_1)|$ (using Proposition \ref{lem:2})\cite{Luo}, \cite{MerrisII}). Hence the proposed graph $G$ is again an LELM for this case also.      
\end{enumerate}
This completes the proof of Case 3 and also the proof of the theorem.
\end{proof}

\section{Main results: globally optimal graphs}\label{sec:GlobOptGraphs}
In this section we obtain sufficient conditions for the union of complete graphs to be a 
global minimizer of the largest eigenvalue. The first main result  of this section (Theorem~\ref{thm:TwoCompLEGM})
is a slight improvement (though claimed and proved on the \emph{complement} graph using different proof
techniques) to Proposition~\ref{lem:BipartiteACM}. The second main result of this section
(Theorem~\ref{thm:MultiCompLEGM}) is a generalization to the case of more than two components and
also gives the first one as a corollary, except the case of equality within the sufficient condition,
Equation~\eqref{eq:TwoCompLEGM}.

\begin{theorem}\label{thm:TwoCompLEGM}
Consider graph $G=(V,E)$ of $n$ number of vertices and $m$ number of edges consisting of two
complete components $K_{\ell}$ and $K_{n-\ell}$, i.e. $m=|E(G)|=\Comb{\ell}{2}+\Comb{n-\ell}{2}$.
Let without loss of generality $\ell \leqslant \frac{n}{2}$. Assume
\begin{fleqn}
\begin{equation}\label{eq:TwoCompLEGM}
\hspace*{1cm} \ell-\frac{2\ell^2}{n}\leqslant 1.
\end{equation}
\end{fleqn}
Then the graph $G=K_\ell\bigcup K_{n-\ell}$ is a Largest Eigenvalue Global Minimizer (LEGM). 
\end{theorem}
\begin{proof}: This proof involves two cases depending on whether the inequality $\displaystyle \ell-\frac{2\ell^2}{n} \leqslant 1$ is strict (Case 1) or holds with equality (Case 2).\\
\noindent \textbf{Case~1: $\displaystyle \ell-\frac{2\ell^2}{n} < 1$}. \\
First notice that when $\displaystyle \ell=\frac{n}{2}$, we get $\displaystyle \ell-\frac{2\ell^2}{n}=0$ and $\displaystyle \ell<\frac{n}{2}$ is same as $\displaystyle 0<\ell-\frac{2\ell^2}{n}$. Hence, the assumption in the theorem gives $\displaystyle 0\leqslant \ell-\frac{2\ell^2}{n}\leqslant 1$. In order to prove the theorem, we obtain the average degree of the graph.\\
Average degree ($d_{avg}$) of the graph $G=K_\ell\bigcup K_{n-\ell}$: 

\begin{equation*}
\begin{split}
    d_{avg} &=\frac{2m}{n}=\frac{2}{n}\left\{\frac{\ell^2-\ell}{2}+\frac{(n-\ell)^2-(n-\ell)}{2}\right\},\\
          &=\frac{1}{n}\left\{2\ell^2+n^2-2n\ell-n\right\},\\
          &=n-\ell-1 -\left(\ell-\frac{2\ell^2}{n}\right).
\end{split}
\end{equation*}

We use that the maximum degree of the graph, $\Delta \geqslant d_{avg}$. In fact, we also use that $\Delta$ should be an integer which implies $\Delta \geqslant \lceil d_{avg}\rceil$. If  $0\leqslant \ell-\frac{2\ell^2}{n}<1$, then the maximum degree, $\Delta \geqslant n-\ell-1$. \\
Using Proposition~\ref{lem:2}a), for any graph that has as many edges as $m$, we get $\lambda_1(G)\geqslant \Delta +1$ and thus $\lambda_1(G) \geqslant n-\ell$ for any graph having as many edges as $G=K_{\ell}\bigcup K_{n-\ell}$. \\
For the proposed graph $G$, the largest eigenvalue of the graph, 
$\lambda_1(G)=\max\{\ell, n-\ell\}=n-\ell$. \\
Hence, the proposed graph $G$ of theorem $K_\ell\bigcup K_{n-\ell}$ is an LEGM. \\
\noindent \textbf{Case~2: $\displaystyle \ell-\frac{2\ell^2}{n} = 1$}. \\
$\displaystyle \ell-\frac{2\ell^2}{n} = 1$ $\implies$ $2\ell^2-n\ell+n=0$\\
whose roots are: $\ell=\displaystyle \frac{n \pm \sqrt{n^2-8n}}{4}$. \\
Notice that for $\ell$ to be an integer the discriminant $n^2-8n$ needs to be a perfect square, i.e. $n^2-8n=p^2$, where $p \in \mathbb{Z}^+$. \\
It is easy to verify that a non-negative integer solution $n$ exists only for $n=9$ in which case $\ell=3$. \\
For this case, i.e. $K_3\bigcup K_6$, we have $\Comb{3}{2}+\Comb{6}{2}=3+15=18$ edges, and $\lambda_1(G)=6$. For this case, through a brute force exhaustive search for $18$ edges, we conclude that $K_3\bigcup K_6$ is an LEGM. (see also Example~\ref{ex:LELM_Equalto_Circulant}).\\
This completes the proof of Theorem~\ref{thm:TwoCompLEGM}.
\end{proof}

We now generalize Theorem~\ref{thm:TwoCompLEGM} and Proposition~\ref{lem:completeComp} to $p$, with $p>2$, components. 

\begin{theorem}\label{thm:MultiCompLEGM}: Consider graph $G$ of $n$ number of vertices and $m$ number of edges consisting of $p$ complete components $K_{n_1}, K_{n_2}, \hdots  K_{n_p}$ such that $\sum_{i=1}^p n_i=n$ and $m=|E(G)|=\sum_{i=1}^p \Comb{n_i}{2}$. Let without loss of generality $n_1\geqslant n_2\geqslant \hdots \geqslant n_p$. Assume \begin{fleqn}
\begin{equation}\label{eq:MultiCompLEGM}
 \displaystyle n_1-\frac{n_1^2+n_2^2+\hdots+n_p^2}{n}<1.   
\end{equation}
\end{fleqn} 
Then the graph $G=K_{n_1}\bigcup K_{n_2}\hdots\bigcup K_{n_p}$ is a Largest Eigenvalue Global Minimizer (LEGM). 
\end{theorem}
\begin{proof}: For the graph $G=K_{n_1}\bigcup K_{n_2}\hdots\bigcup K_{n_p}$, 
first notice that $n_1-\frac{n_1^2+n_2^2+\hdots+n_p^2}{n} \geqslant 0$.\\
This is because 
\[
\Scale[1.0]{\frac{n\times n_1-(n_1^2+n_2^2+\hdots+n_p^2)}{n}=\frac{n_2(n_1-n_2)+n_3(n_1-n_3)+\hdots + n_p(n_1-n_p)}{n}}
\]
and thus only when $\frac{n}{p} \in \mathbb{Z}$ (and hence
$\frac{n}{p}=n_1 = n_2 = \hdots = n_p$, we have $n_1-\frac{n_1^2+n_2^2+\hdots+n_p^2}{n}=0$.
For any other value of $n$ and of $n_{i}$,
we have  $0<n_1-\frac{n_1^2+n_2^2+\hdots+n_p^2}{n}$ and thus  ${0\leqslant n_1-\frac{n_1^2+n_2^2+\hdots+n_p^2}{n}}$ in general.\\
The average degree ($d_{avg}$) of the graph: 
\begin{equation*}
\begin{split}
    d_{avg} &=\frac{2m}{n}=2\sum_{i=1}^p\frac{\Comb{n_i}{2}}{n}, \\
          &=\sum_{i=1}^p\frac{n_i^2-n_i}{n}=\frac{n_1^2+n_2^2+\hdots +n_p^2-n}{n},\\
          &=n_1-1 -(n_1- \frac{n_1^2+n_2^2 \hdots +n_p^2}{n} ).
\end{split}
\end{equation*}
We next use that the maximum degree of the graph, $\Delta \geqslant d_{avg}$. We also know that $\Delta$ should be an integer which implies $\Delta \geqslant \lceil d_{avg}\rceil$. 
If ${0\leqslant n_1-\frac{n_1^2+n_2^2 \hdots +n_p^2}{n}<1}$, then the maximum degree $\Delta \geqslant n_1-1$. \\
Using Proposition~\ref{lem:2}a), for any graph that has as many edges as $m$, we get: $\lambda_1(G)\geqslant \Delta+1$ $\implies$ $\lambda_1(G)\geqslant n_1$.\\ 
Finally, it remains to show that the proposed graph $G$ satisfies $\lambda_1(G)=n_1$. Since, $n_1\geqslant n_i$ and $\lambda_1(K_{n_i})=n_i $, we conclude that graph proposed is LEGM. This completes the proof of Theorem~\ref{thm:MultiCompLEGM}.
\end{proof}

\begin{remark} \label{rem:sufficiency:two:components}
Theorem~\ref{thm:TwoCompLEGM} (and Theorem~\ref{thm:MultiCompLEGM}, for the number of components greater than two case) establish that when two components are of `almost similar sizes' or the largest component is `much larger than the smallest', we get a Largest Eigenvalue Global Minimizer graph.
Both sufficient conditions, equations (\ref{eq:TwoCompLEGM}) and (\ref{eq:MultiCompLEGM}) are to be viewed as a relaxation on the condition `$\alpha$ is factor of $n$' in Proposition \ref{lem:completeComp}.
This is elaborated as follows. From Proposition~\ref{lem:completeComp}, it is clear that for any integer $n_1$, when we have $\bigcup_{i=1}^p K_{n_{1}}$, then
this graph is an LEGM.
Intuitively, by addition of a `sufficiently small' component $K_{n_{p+1}}$, i.e.
when $0<n_{p+1} ~{\large{ \mathlarger{\ll}}} ~n_1$, then LEGM would continue to hold.
On the other hand, when $n_p$ is `slightly smaller' than $n_1$, then too LEGM would continue to hold.
In other words, not just when all components are of same size but when the components
are `quite homogeneous' or some of them are relatively `negligibly small', then also LEGM property continues to
hold: in that sense the sufficient condition $\displaystyle n_1-\frac{n_1^2+n_2^2+\hdots+n_p^2}{n}<1$ is a relaxation of the condition `$\alpha$ is factor of $n$' of Proposition~\ref{lem:completeComp}. 
\end{remark}
 
\begin{theorem}
Consider the algorithm below that takes $n$ (number of vertices) and $m_{desired}$ (desired number of edges) as an input.
Suppose the algorithm terminates with $m_{actual}=m_{desired}$, then the constructed graph is LELM.
If the sufficient condition of Theorem~\ref{thm:MultiCompLEGM} is met, then
this proposed graph is also LEGM.
Within the class of graphs which are LELM, this procedure gives the least $\lambda_1$.
\end{theorem}
\begin{proof}
The claims in the theorem are straight forward and hence we summarize and
dwell on only the key arguments.
The algorithm constructs components: largest first and then smaller, etc.
until all vertices are
used up and the maximum number of edges (up to $m_{desired}$) are accommodated.
\begin{itemize}
\item By construction, the obtained graph is clearly LELM. 
\item Within the `the while loop', the condition $x_i\leqslant n_i^{rem}$ ensures that the new components do not exceed the remaining number of vertices. 
\item Equation (\ref{eq:algL2}) ensures that $\ell_i$ is as large as possible for a given component size $x_i$.
\item Equation (\ref{eq:algL1}) ensures that the $\ell_i$ components, each of $x_i$ vertices, do not exceed the remaining number of vertices. 
\end{itemize}
Thus, the construction procedure attempts to accommodates the desired number of edges with as small size components of complete graphs $K_{x_i}$ as possible and hence is an LELM with least $\lambda_1=x_1$. 
\end{proof}

\begin{algorithm} \label{alg:LELM:construction:ajay}
\SetAlgoLined
\textbf{Input}: {\small Vertices count: $n$, number of edges desired: $m_{desired}$.} \\
\textbf{Output}: {\small Number of components $p$, their sizes and the number of edges actually accommodated $m_{actual}$.}

\textbf{Initialize} $i=1$, the number of vertices in the graph $n_{1}^{rem}:=n$ and the desired number of edges to be accommodated in the graph $E_1^{rem}:=E_{d_{1}}=m_{desired}$. \\
 \While{$E^{rem}_{i+1} > 0$ \emph{~and~} $n^{rem}_{i+1} > 0$,~}{
 Get $x_i$, $\ell_i$ by the following minimization:
  $\arg \min x_i \in \mathbb{Z}^+$, $x_i \leqslant n^{rem}_{i}$ \\
  such that there exists $l_i \in \mathbb{Z}^+$ satisfying (\ref{eq:algL1}) \& (\ref{eq:algL2}):\\
  \begin{fleqn}[\parindent]
  \begin{equation}\label{eq:algL1}
      \hspace*{2mm}  \ell_{i}\leqslant \frac{n^{rem}_{i}}{x_{i}} 
  \end{equation}
  \vspace*{-7mm}
  \end{fleqn}
  \begin{fleqn}
  \begin{equation}\label{eq:algL2}
      \hspace*{4mm} \ell_i~\Comb{ x_i }{2}\leqslant E^{rem}_{i} <  (\ell_i+1)~\Comb{x_{i}}{2}
  \end{equation}
  \end{fleqn}
  $E^{rem}_{i+1}:=E^{rem}_{i}-\ell_i~\Comb{x_i}{2}$ and $n^{rem}_{i+1}:=n^{rem}_{i}-\ell_ix_i$\;
 }
\KwResult{Suppose at $i=s$, one or both of the conditions, $E^{rem}_{i+1} > 0$  or  $n^{rem}_{i+1} > 0$ gets violated,
  then the following are defined as the output. 
\begin{enumerate}
    \item The actual number of edges the constructed graph accommodates, $m_{actual}:=\sum_{i=1}^{s} \ell_i~\Comb{x_i}{2}=|E(G)|$.
    \item The number of components of the graph, $p=\sum_{i=1}^{s} \ell_i$.
    \item The LELM graph $G:=\underset{\ell_1 \mbox{times}}{K_{x_1}\bigcup} K_{x_2} \bigcup\cdots \bigcup\underset{\ell_{s} \mbox{times}}{K_{x_{s}}} $, and $\lambda_1(G) = x_1$
\end{enumerate}}
 \caption{Edges inclusion in the graph}
\end{algorithm} 

Note that when the algorithm terminates, but with $m_{desired}>m_{actual}$, then
the difference $m_{desired}-m_{actual}\leqslant n-2$: this can be seen easily and is
hence not pursued. Obtaining a better upper bound is worth pursuing further.

\begin{remark} \label{rem:sufficient:interpretation:global}
 Equations~(\ref{eq:algL1}) and (\ref{eq:algL2}) within Algorithm~\ref{alg:LELM:construction:ajay} are to be understood as follows.
It is understandable that to have $\lambda_1$ low, the complete components need to be individually of a small sizes. 
This is achieved by taking the min $x_i$ satisfying equations $(\ref{eq:algL1})$ and $(\ref{eq:algL2})$.
The condition $x_i\leqslant n^{rem}_{i}$ is about how many vertices are available
for the next graph construction.
For each $x_i$, the number of components $K_{x_i}$ is $\ell_i$.
Each component of size $x_i$ accommodates $\Comb{x_i}{2}$ edges and we try to have as many such components of size $x_i$ as possible given the total remaining number of vertices $n^{rem}_{i}$, this is captured by $\ell_i\leqslant\frac{n_{i}^{rem}}{x_i}$.
Finally, given a size $x_i$, it is required that the number of components $\ell_i$ of that size should be as large as possible to accommodate the desired (or yet to be accommodated) number of edges: this is ensured by Equation~(\ref{eq:algL2}).
Loosely speaking, increasing $x_i$ helps in accommodating more edges, at the cost of a larger $\lambda_1$ and
less number of components $\ell_i$. On the other hand, smaller $x_i$ aids in decreasing $\lambda_1$, but would perhaps be unable to accommodate enough edges.
\end{remark}

\section{Circulant matrices}\label{sec:CirculantMatrices}

Propositions~\ref{lem:1} and~\ref{lem:2} are about relations between
the Laplacian eigenvalues for a graph and its complement, and about
the max degree $\Delta$ providing a lower bound for the max eigenvalue.
In particular, the lower bound $\Delta + 1$ is tight for the case
when the graph contains a star node, i.e. the domination number (see
Footnote~\ref{foot:domination:number}) is 1. This naturally suggests that a
relatively equitable distribution of edges that keeps the max-degree
$\Delta$ low helps in keeping the maximum eigenvalue $\lambda_1$
also low.

Circulant matrices are such matrices: they are regular and contain a symmetry
that indeed makes them LEGM for certain cases; we pursue this link
in this section.

A matrix $C\in\Rnn$ is called circulant if each entry $c_{i,j}$, the entry in $i$-th row
and $j$-th column satisfies:
$c_{i,j} = c_{i+k,j+k}$, where the indices are considered to be modulo-$n$ and - for
this reason, and just for this sentence - indices $i,j$ vary from $0$ to $n-1$. It is well-known (see \cite{Philip-Davis})
that the set of circulant matrices form
an $n$-dimensional subspace of $\Rnn$, and the entries of only the first row of $C$ need to
be specified for specifying $C$.
A circulant graph is one whose Laplacian is a
circulant matrix, after a permutation/re-ordering of the nodes, if needed.
Define the matrix $J \in \mathbb{R}^{n \times n}$ such that $J_{ij}=1$  for all $i, j \in \{1, 2,\cdots, n\}$.
Notice that $nI-J$ is a circulant matrix with generating row as $[n-1, -1, -1, \cdots, -1]$.
The Laplacian of this circulant graph is same as the Laplacian of $K_n$,  i.e. $nI-J$.
This means that if $G^m$ is a circulant graph, then so is its complement $G^c$.

We pursue further with Problem~\ref{prob:DFT} and note that the DFT of the first row
of a circulant matrix $C$ are exactly the eigenvalues of $C$.
Given integers $n$ and $m$, the number of vertices and edges, due to the implicit
regularity of a circulant graph, $2\times m$ has to be divisible by $n$ for a circulant
graph $G(V,E)$ to exist such that $|V| = n$ and $|E|=m$.

Below is our first result in this context. We then come up with examples in the
following section.

\begin{theorem} \label{thm:circulant:LEGM}
Consider positive integers $n$ and $m$ satisfying the relation that $n$ is a factor of $2m$, and $\frac{2m}{n}+1$ is a factor of $n$.
Then, the following hold.
\begin{enumerate}
\item There exists a circulant graph $G_c^m$ having $n$ vertices and $m$ edges.
\item $G_c^m$ is an LEGM.
\item The first row of the adjacency matrix of $G_c^m$ solves Problem~\ref{prob:DFT}.
\item $G_c^c$, the complement of $G_c^m$,  is also a circulant graph and has the highest
  algebraic connectivity, i.e. $G_c^c$ is an ACM.
\end{enumerate}
\end{theorem}

\begin{proof}
Notice that the condition on $m$ is just that one can construct
$G:=\underset{\ell~\mbox{times}}{\bigcup K_i}$, with $i:=\frac{2m+n}{n}$ and $\ell:=\frac{n}{i}$.
By a careful renumbering of the vertices in $G$,
it is possible to obtain a circulant graph $G_c$. Note that renumbering of vertices
is merely premultiplying and postmultiplying the Laplacian $L$ by permutation matrices
$P$ and $P^T$, a unitary similarity transformation, does not change
the eigenvalues of $L$.
\end{proof}

Of course, the condition specified in the theorem is only a sufficient condition for
a circulant matrix to be an LEGM.
Examples~\ref{ex:LELM_Equalto_Circulant} and~\ref{ex:LELM_Betterthan_Circulant}
are included in the next section:
The former (i.e. Example~\ref{ex:LELM_Equalto_Circulant}) is
when the sufficient condition of Theorem~\ref{thm:TwoCompLEGM}
is met with equality, and, and the resulting union of complete components
is an LEGM. Further, this case also admits a circulant matrix, which also is an LEGM,
though it is not a union of complete components.
The latter (Example~\ref{ex:LELM_Betterthan_Circulant}) is a circulant matrix that has
$\lambda_1$ significantly higher than the corresponding LELM constructed
for $n=9$ and $m=18$.

\section{Examples} \label{sec:examples}

In this section we consider some examples.  Table~\ref{table:nmvalues:LEGM:LELM} contains
many typical values of $n$ and $m$ (the total number of vertices and edges) and also lists which are LEGM (in addition to being LELM). Some more examples
are elaborated here.

\begin{table}[!h] 
\centering
\caption{$\lambda_1$ for complete components $K_i$ having $n$ vertices and $m$ edges} \label{table:nmvalues:LEGM:LELM}
    \begin{tabular}{||c|c|c|c|c||}
    \hline
    $n$ & $m$ &	$\bigcup K_i$	& $\lambda_1$ & LEGM/LELM\\
    \hline
    9& 10&	4,  3,  2	& 4 &LEGM\\
    9& 12&	4,  4 &4 &LEGM\\
    10& 16&	5,  4  &5  &LEGM\\
    10& 20&	5,  5	 &5 &LEGM\\
    15& 34&	7,  5,  3	&7 &LELM\\
    %
    20& 22&	4, 4, 4, 2, 2, 2, 2&	4 &LEGM\\
    20& 50&	8, 6,  4,  2&	8 &LELM\\
    25& 66&	8, 8, 4, 3, 2	&8 &LELM\\
    25& 132&	12, 12&	12 &LEGM\\
    30& 235&	22, 2, 2, 2, 2	&22 &LELM\\
    32& 136& 10, 10, 10, 2 &10 &LEGM\\
    \hline
    \end{tabular}
    \label{tab:EdgeAcmmdtKn}
\end{table}

\begin{example}\label{ex:LELM_Equalto_Circulant}
In this example, the sufficient condition Inequality~(\ref{eq:TwoCompLEGM}) is satisfied with an equality
and is not captured by Proposition~\ref{lem:BipartiteACM}, but handled by Theorem~\ref{thm:TwoCompLEGM}.
Suppose the number of vertices,~$n=9$ and the number of edges,~$m=18$.
\begin{itemize}
\item The LELM graph is $K_6\bigcup K_3$ with $\lambda_1(K_6\bigcup K_3)=6$, and by a simple exhaustive
brute-force search, this also is an LEGM.
\item Further, the circulant graph $G_c$ with degree 4, represented
by the circulant adjacency matrix having its first row as
$\left[0~ 0~ 1~ 0~ 1~ 1~ 0~ 1~ 0\right]$ also has $\lambda_1(G_c)=6$. 
Thus $K_6\bigcup K_3$ is not the unique LEGM and
the circulant graph $G_c$ has the same $\lambda_1$ value and is an LEGM too.
\end{itemize}
\end{example}

\begin{example}\label{ex:LELM_Betterthan_Circulant}
Consider again the case when vertex/edge counts are $n=9$ and $m=18$,
and we look for
a circulant graph that {\em maximizes} the largest eigenvalue.  As noted in the previous
example, the LELM graph $K_6\bigcup K_3$ gives $\lambda_1(LELM)=6$.
A different  circulant graph $G_c$, obtained from the circulant adjacency matrix having its first
row as $\left[0~0~ 0~1~ 1~ 1~ 1~ 0~ 0\right]$, has $\lambda_1(G_c)=6.88$.
\end{example}

\begin{example}\label{ex:LELM_worstthan_Circulant}
Consider the case when the vertex/edge counts are~$n=7$ and~$m=7$.
\begin{itemize}
   \item LELM graph is $K_4\bigcup K_2\bigcup K_1$ has $\lambda_1(LELM)=4$.
   \item The circulant graph $C_7$, the cycle graph on 7 nodes, represented by circulant adjacency matrix generated by $\left[0~1~0~0~0~0~1\right]$ has $\lambda_1(C_7)=3.802$.
\end{itemize}
\end{example}

\begin{example}\label{ex:LELM_worstthan_Circulant2}
Consider the case when the vertex/edge counts are $n=24$, $m=168$.
In this case, the sufficient condition of Theorem~\ref{thm:TwoCompLEGM} is violated relatively quite severely. The LELM graph generated by our algorithm is $K_{18} \bigcup K_6$.
This is a case where the two components are too heterogeneous, and the LELM graph is not LEGM.
The circulant graph $G_c$ obtained by the circulant adjacency matrix having its
first row as $\left[
0~0~0~1~1~0~1~1~1~0~1~1~0~1~1~0~1~1~1~0~1~1~0~0 
\right]$ has $\lambda_1(G_c)=17$. 
\end{example}

\section{Concluding remarks}\label{sec:conclusion}
In this paper we showed how the graphs comprised of two or more complete components 
locally minimize the Laplacian's largest eigenvalue (LELM graphs): Theorems~\ref{thm:MultiCompLELM}. This was achieved by a meticulous case by case analysis of various possibilities: see Figures~\ref{fig:Case1:1a} to~\ref{fig:Case3:3r}, and Lemmas~\ref{lem:AddEdge}-\ref{lem:ReconnectEdge_Sdmax}.
Further, if the components sizes are either `quite homogeneous' or some of them are relatively
`negligibly small' (as elaborated in Remark~\ref{rem:sufficient:interpretation:global}, which
interpreted Equations~(\ref{eq:TwoCompLEGM}) and (\ref{eq:MultiCompLEGM}) of Theorems~\ref{thm:TwoCompLEGM}
and~\ref{thm:MultiCompLEGM}), then this graph is not just local, but
also a global minimizer of the largest eigenvalue for that many vertices and edges. This thus extends existing results in different and appropriate ways: Propositions~\ref{lem:1}, \ref{lem:2} and~\ref{lem:DominatingNum}. 
We also proposed
an algorithm to construct such a locally/globally optimum
graph (Algorithm~\ref{alg:LELM:construction:ajay}).

We also related our results to the well-studied class of graphs called circulant graphs: the significance
being that due to the symmetry and fairly uniform distribution
of edges across nodes within such graphs, they appear like the opposite of graphs that
have a `star node', and hence are potential candidates for minimization of the largest eigenvalue.
The link between circulant graphs/matrices and the Discrete Fourier Transform is well-known, and
the central problem considered in this paper thus translates to minimization of the maximum
magnitude across all nonzero frequencies in a periodic discrete time signal (see Problem~\ref{prob:DFT}
and Remark~\ref{rem:DFT:Circulant}).

Our method of maximizing the algebraic connectivity of a graph crucially uses Proposition~\ref{lem:1}.
Thus all our results pertaining to minimization of the largest eigenvalue easily 
translate to the maximizing of the algebraic connectivity by noting that the main graph $G^m$ and its complement $G^c$ (on $n$-vertices) have the edge counts adding up to ${}^n C_2$, i.e. $\frac{n(n-1)}{2}$.

\bibliographystyle{IEEEtran}
\bibliography{refs}


\end{document}